\documentclass[11pt, a4paper]{article}
\def\qed{\ifmmode\square\else\nolinebreak\hfill$\diamondsuit$\fi\par\vskip12pt}

\newcommand{\be}{\begin{equation}}
\newcommand{\ee}{\end{equation}}
\newcommand{\DS}{\displaystyle}

\newcommand{\beq}{\begin{eqnarray}}
\newcommand{\eeq}{\end{eqnarray}}
\newcommand{\nbeq}{\begin{eqnarray*}}
\newcommand{\neeq}{\end{eqnarray*}}

\newcommand{\AmS}{{\protect\the\textfont2
  A\kern-.1667em\lower.5ex\hbox{M}\kern-.125emS}}

\def\be{\begin{equation}}
\def\ee{\end{equation}}

\begin{document}
\title{\bf Number of complete $N$-ary subtrees on Galton-Watson family trees}
\author{\small
GEORGE P. YANEV \hspace{\fill} gyanev@cas.usf.edu \\
\small \hspace{-2.3cm} Department of Mathematics, University of
South Florida,
Tampa, FL 33620, USA\\
\small LJUBEN MUTAFCHIEV  \hspace{\fill} ljuben@aubg.bg \\
\small American University in Bulgaria, 2700
 Blagoevgrad,
 Bulgaria and  Institute of Mathematics and \\ \small \hspace{-7.2cm}  Informatics of the Bulgarian Academy of Sciences}
\date{\empty}

\maketitle

\begin{abstract}

We associate with a Bienaym\'{e}-Galton-Watson branching process a
family tree rooted at the ancestor. For a positive integer $N$,
define a complete $N$-ary tree to be the family tree of a
deterministic branching process with offspring generating function
$s^N$. We study the random variables $V_{N,n}$ and $V_N$ counting
the number of disjoint complete $N$-ary subtrees, rooted at the
ancestor, and having height $n$ and $\infty$, respectively.
Dekking (1991) and Pakes and Dekking (1991) find recursive
relations for $P(V_{N,n}>0)$ and $P(V_N>0)$ involving the
offspring probability generation function (pgf) and its
derivatives. We extend their results determining the probability
distributions of $V_{N,n}$ and $V_N$. It turns out that they can
be expressed in terms of the offspring pgf, its derivatives, and
the above probabilities. We show how the general results simplify
in case of fractional linear, geometric, Poisson, and one-or-many
offspring laws.
\end{abstract}

\hspace{0.5cm}{\bf Keywords:} branching process - family tree -
binary tree - $N$-ary tree

\hspace{0.5cm}{\bf AMS 2000 Subject Classification:} Primary:
60J80. Secondary: 05C05.

\section{Introduction and main results}

Consider the family tree associated with a
Bienaym\'{e}-Galton-Watson process with the following simple
reproduction rules. At generation zero, the process starts with
single ancestor called root of the tree. Then each individual in
the population has, independently of the others, a random number
$X$ of children distributed according to the offspring
distribution with probability generating function (pgf)
\[
f(s)=\sum_{k=0}^\infty p_k s^k,
\]
satisfying $f(1)=1$. Further on we adopt the well-known
construction of a family tree generated by a simple branching
process where the individuals are the nodes and the parent-child
relations define the arcs of the tree in the following manner, see
e.g. Harris (1963), Ch.7. Let the $i$th child of the ancestor be
$(i)$ and in general $(i_1i_2\ldots i_{k-1}i_k)$ denotes the
$i_k$th child of $(i_1i_2\ldots i_{k-1})$. Then, a directed arc is
assumed to emanate from $(i_1i_2\ldots i_{k-1})$ to $(i_1i_2\ldots
i_{k-1}i_k)$. Since, in our case, the children appear
simultaneously, we suppose that the ordering is performed by a
chance device independently of the evolution in the process. This
scheme produces family trees (also called rooted ordered trees) in
which the nodes of height (also known as depth) $n$ $(n\geq 0)$
have labels $(i_1i_2\ldots i_n)$, with the ancestor (root) having
height 0. The height of a subtree equals the maximum height of its
nodes.

For fixed integer $N\ge 1$, define a complete infinite $N$-ary
tree to be the family tree of a deterministic branching process
with offspring pgf $f(s)=s^N$. Further on we will consider rooted
subtrees of a family tree. Two such subtrees are called disjoint
if they do not have a common node different from the root. These
kinds of trees appear, for example, in some computer algorithms;
for more details see Knuth (1997).

Let $\{Z_n: n\geq 1; \ Z_0=1\}$ denote the generation size
process, and let $T_N-1$ be the height of a complete $N$-ary
subtree rooted in the ancestor; $T_N=0$ if $Z_1<N$. Notice that
$T_1$ is the extinction time of $\{Z_n\}$. The study of the
probability $\tau_N=\lim_{n\to \infty}P(T_N>n)$ that a
Bienaym\'{e}-Galton-Watson tree contains an infinite complete
$N$-ary subtree was initiated by Dekking (1991) who considered
complete binary ($N=2$) subtrees. The general $(N\geq 2)$ case was
subsequently investigated in detail by Pakes and Dekking (1991).
In particular, they encountered the following phenomenon: if
$N\geq 2$, then there is a critical value $m^c_{N}$ for the
offspring mean $m=f'(1)$ such that $\tau_N=0$ if $m<m^c_N$ and
$\tau_N>0$ if $m\ge m^c_{N}$. This is qualitatively different from
what happens for $N=1$ where the probability for non-extinction
$\tau_1=0$ if $m=m^c_1=1$, except for the trivial case where
$f(s)=s$. Our work is motivated by the results of Pakes and
Dekking (1991).

We introduce the random variable $V_N $ to be the number of
disjoint complete $N$-ary subtrees with infinite height, rooted at
the ancestor of a Bienaym\'{e}-Galton-Watson family tree. Clearly
$\tau_N=P(V_N>0)$.  As usual, we assume for the offspring
distribution $\{p_k\}_{k=0}^\infty$ that $p_k<1$ for all $k$ and
$p_k>0$ for some $k>N$. Let $\cal N$ be the set of all positive
integers and denote for $x,y\geq 0$ and any $j=0,1,\ldots$
\[
G_N(x,y;j)=\sum_{k=jN}^{jN+N-1} \frac{x^{k}}{k!}f^{(k)}(y).
\]
 Pakes and Dekking (1991) showed that $P(V_N=0)=1-\tau_N$, where
$1-\tau_N$ is the smallest solution in $[0,1]$ of the equation \be
\label{first_eq} x=G_{N}(1-x,x;0). \ee

Our goal is to study the distribution of $V_N$. As the following
result shows, the probability mass function (pmf) of $V_N$ can be
obtained using the Taylor expansion of $f(1)$ about the point
$1-\tau_N$.

{\bf Theorem 1} If $N\in \cal N$ then for any $j=0,1,...$ \be
\label{main} P\left(V_N=j)=G_N(\tau_N, 1-\tau_N;j\right)
 \ee and $P(V_N=0)=1-\tau_N$ is the smallest solution
in $[0,1]$ of (\ref{first_eq}).

\emph{Remark}  (i) If $N=1$, then obviously $P(V_1=0)=1-\tau_1=q$
is the extinction probability of the Galton-Watson process. Now,
(\ref{main}) becomes
\[
P(V_1=j)=\frac{(1-q)^j}{j!}f^{(j)}(q), \qquad j=0,1,\ldots,
\]
which in turn implies that $E(s^{V_1})=f(q+(1-q)s)$. This identity
follows directly observing that the number of distinct infinite
unary trees is equal to the number of first generation nodes
having infinite line of descent.\footnote{ The authors are
indebted to the referee who pointed out this argument. It implies
immediately the result of Theorem 1 for unary trees.}

(ii) Also note that a sufficient condition for $P(V_N=0)<1$ is
given in Pakes and Dekking (1991), Theorem 3. In particular, they
show that $P(V_N=0)<1$ $(N\geq 2)$ if
$$
2N\sum_{j\ge N}\frac{p_j}{j+1}\le (1-\sum_{j=0}^{N-1} p_j)^2.
$$

The number of complete $N$-ary subtrees is a measure for the rate
of growth (or fertility) of the branching process. In fact, as was
pointed out in Dekking (1991), if $P(V_2>0)>0$ then we can say
that the branching process grows faster than binary splitting. In
the study of the tree structure of branching processes, an
important role is played by the process' total progeny. Denote by
$\nu_n$ the number of individuals who existed in the first $n+1$
generations, i.e., $\nu_n=1+Z_1+\ldots +Z_n$, $n=1,2,\ldots$.
Obviously, $\nu_n$ equals the total number of nodes having height
less than or equal to $n$. Let us also define the random variable
$V_{N,\ n}$ to be the number of disjoint complete $N$-ary subtrees
of height at least $n$ rooted at the ancestor of a
Bienaym\'{e}-Galton-Watson family tree. Let \be \label{psiphi}
\psi_{N,n}(s)=E(s^{\nu_n}; V_{N,n}>0) \quad \mbox{and} \quad
\phi_{N,n}(s)=E(s^{\nu_n}; V_{N,n}=0). \ee
 The following result
presents a recursive relation for the joint distribution of
$V_{N,n}$ and $\nu_n$.

{\bf Theorem 2}\ If $N\in \cal N$ then for $\mid s\mid\leq 1$ and
any $j=0,1,...$ \be \label{recurrence} E\left(
s^{\nu_{n+1}};V_{N,n+1}=j\right)=s
G_N\left(\psi_{N,n}(s),\phi_{N,n}(s);j\right). \ee

Notice that, if $N=1$ and $j=0$, then the above recurrence reduces
to the well-known $E\left( s^{\nu_{n+1}}; Z_{n+1}=0\right)=s
f\left(E\left( s^{\nu_{n}};Z_{n}=0\right)\right)$, see e.g.
Kolchin (1986), p. 120.


 Applications of complete $N$-ary trees can be found in the
analysis of algorithms, see Knuth (1997). Problems of this nature
appear also in percolation theory. For instance, Pakes and Dekking
(1991) point out a relationship between the model of $N$-ary
complete and infinite subtrees and a construction employed by
Chayes et al. (1988) in their study of Mandelbrot's percolation
processes. The existence of $N$-ary subtrees is also used by
Pemantle (1988) in introducing the concept of a $N$-infinite
branching process. Let us also mention potential connections with
problems of percolation of binary words on the nodes of locally
finite graphs with countably infinite node-sets, see Benjamini and
Kesten (1995).

We organize our paper as follows. In Section 2 we prove the main
results. Sections 3-5 contain some illustrations. In Section 3 we
consider the family tree generated by the fractional linear $f(s)$
as well as the special case of geometric offspring. In the latter
case, $V_N$ itself follows a geometric distribution. It turns out
that in the Poisson offspring case, given in Section 4, the pmf of
$V_N$ can be expressed in terms of certain Poisson probabilities.
Note that the critical values $m^c_N \ (N\geq 2)$ in the Poisson
case are less than those in the geometric one. Finally, in
Section~5 we consider the one-or-many (i.e., concentrated on two
points only) offspring distribution. In this case $V_N$ has a pmf
given in terms of binomial probabilities.

\section{Proofs of the Theorems}

{\bf Proof of Theorem 1}\ \
Let us consider $P(V_N=j)$
where $j=1,2,\ldots$ Recall that the random variable $V_{N,\ n}$
equals the number of disjoint complete $N$-ary subtrees of height
$n$ rooted at the ancestor of a Bienaym\'{e}-Galton-Watson family
tree. First, we will find the pmf of $V_{N,\ n+1}$ using the total
probability formula. Indeed, to have $j$ disjoint complete $N$-ary
subtrees rooted at the ancestor node there must be $jN+k \ (k\geq
0)$ nodes in the first generation.
 Each of these nodes can be considered as an
ancestor of a family tree rooted at the first generation. Consider
the event $A_N(l)= \{jN+l$ of the $Z_1$ first generation nodes are
ancestors of at least one complete $N$-ary tree of height $n\}$,
where $l=0,1,\ldots, \min\{k, N-1\}$. If $Z_1=jN+k$ then for fixed
$l$ the event $A_N(l)$ has conditional probability
\[
P(A_N(l)|Z_1=jN+k)={jN+k \choose jN+l} (\tau_{N,\
n})^{jN+l}(1-\tau_{N,\ n})^{k-l} \quad (0\leq l\leq \min\{k, \
N-1\}),
\]
where $\tau_{N,\ n}=1-P(V_{N,\ n}=0)$ and by convention let
$\tau_{N,\ 0}=1$. We have
\[
P\left(\bigcup_{l=0}^{\min\{k,\ N-1\}}A_N(l)|Z_1=jN+k
\right)=\sum_{l=0}^{\min\{k,\ N-1\}}{jN+k \choose jN+l} (\tau_{N,\
n})^{jN+l}(1-\tau_{N,\ n})^{k-l}.
\]
Applying the total probability formula and changing the order of
summation, we obtain \beq  P(V_{N,\ n+1}=j) & = &
\sum_{k=0}^\infty P(Z_1=jN+k)\
P\left(\bigcup_{l=0}^{\min\{k,\ N-1\}}A_N(l) \ |\  Z_1=jN+k\right) \nonumber \\
    & = &
\sum_{k=0}^\infty p_{jN+k}
    \left\{\sum_{l=0}^{\min\{k,\ N-1\}}
{jN+k \choose jN+l} (\tau_{N,\ n})^{jN+l}(1-\tau_{N,\
n})^{k-l}\right\} \nonumber \\
   & = &
\sum_{l=0}^{N-1}\frac{\tau_{N, n}^{jN+l}}{(jN+l)!}
\sum_{k=l}^\infty
p_{jN+k}(jN+k)(jN+k-1)...(k-l+1)(1-\tau_{N, n})^{k-l} \nonumber \\
    & = &
\sum_{l=0}^{N-1}\frac{\tau_{N,
n}^{jN+l}}{(jN+l)!}f^{(jN+l)}(1-\tau_{N,n}) \nonumber \\
    & = &
    G_N(\tau_{N,n},1-\tau_{N,n};j) \nonumber .
\eeq
 By definition $\tau_{N,\ 0}=1$ and $\tau_{N,\ n}\downarrow \tau_N$ as $n\uparrow \infty$.
 Letting $n\to \infty$, we
 obtain for $j\geq 1$
\[ P(V_N=j) =  \lim_{n\to \infty}P(V_{N,\ n+1}=j)=G_N(\tau_N,
1-\tau_N; j).  \] Let us now consider the case $j=0$. The above
recurrence is true for $n=0$, i.e., $
P(V_{N,1}=0)=G_{N}(1,0;0)=\sum_{k=0}^{N-1}p_k. $ For $n\ge 1$,
using the total probability formula and an argument similar to
that for the case $j\ge 1$, we obtain
\begin{eqnarray}
P(V_{N,n+1}=0) & = & \sum_{l=0}^{N-1}\sum_{k=l}^\infty
p_k{k\choose l}
(\tau_{N,n})^l(1-\tau_{N,n})^{k-l} \nonumber \\
 &  = & \sum_{l=0}^{N-1}\frac{(\tau_{N,n})^l}{l!}
 f^{(l)}(1-\tau_{N,n})\nonumber \\
 & = & G_N(\tau_{N,n},1-\tau_{N,n};0).
 \end{eqnarray}
Computing the derivative of $G_N(x,1-x;0)$, we get a telescoping
sum which after cancelations becomes
$dG_N(x,1-x;0)/dx=(1-x)^{N-1}f^{(N)}(x)/(N-1)! \geq 0$ for $0\le
x\le 1$. Thus, $G_N(x,1-x;0)$ is non-decreasing in $[0,1]$, and
therefore
\[
1-\tau_N=\lim_{n\to \infty}(1-\tau_{N,n+1})=\lim_{n\to
\infty}P(V_{N,n+1}=0) =G_N(\tau_{N},1-\tau_{N};0) \] is the
smallest root in $[0,1]$ of the equation $x=G_N(1-x,x;0)$. The
proof is complete.

Clearly (\ref{main}) implies that $\sum_{j=0}^\infty
P(V_N=j)=\sum_{k=0}^\infty \tau_N^kf^{(k)}(1-\tau_N)/k!=f(1)=1$.

{\bf Proof of Theorem 2} \  Let us introduce the notation
\[
\tau_{N,n}(t)=P(V_{N,n}>0, \nu_n=t), \quad
\gamma_{N,n}(t)=P(V_{N,n}=0, \nu_n=t)=P(\nu_n=t)-\tau_{N,n}(t),
\]
 where $N$, $n$, and $t$ are positive integers. Proceeding as in the proof of Theorem 1, we consider the
event
$$
A_N(l,t)=A_N(l)\bigcap\{\nu_{n+1}=t\},
$$
where $A_N(l)$ is defined in the proof of Theorem 1. For fixed $t$
and $l$ ($0\le l\le \min{(k,N-1)}$), using the fact that all trees
rooted in the first generation grow independently, we compute the
 conditional probability of $A_N(l,t)$ given $Z_1=jN+k$ to
be
$$
P(A_N(l,t)|Z_1=jN+k)={jN+k\choose jN+l}\sum { ^\prime}
\prod_{u=1}^{jN+l}\tau_{N,n}(n_u)
\prod_{v=jN+l+1}^{jN+k}\gamma_{N,n}(n_v),
$$
where the summation in $\sum ^{\prime}$ is over all nonnegative
integers $\{n_i\}_{i=1}^{jN+k}$ such that
$\sum_{i=1}^{jN+k}n_i=t-1$.
 Then, the total probability formula
implies that
\begin{eqnarray}
P(V_{N,n+1}=j,\nu_{n+1}=t) & = & \sum_{k=0}^\infty
P(Z_1=jN+k)\sum_{l=0}^{\min{(k,N-1)}}P(A_N(l,t)\mid Z_1=jN+k) \nonumber \\
& = & \sum_{l=0}^{N-1}\sum_{k=l}^\infty p_{jN+k}{jN+k\choose jN+l}
\sum {^\prime}\prod_{u=1}^{jN+l}\tau_{N,n}(n_u)
\prod_{v=jN+l+1}^{jN+k}\gamma_{N,n}(n_v). \nonumber
\end{eqnarray}
Multiplying both sides of this equality by $s^t$ and summing over
$t$, we get
\begin{eqnarray}
E(s^{\nu_{n+1}}; V_{N,n+1}=j) & = &
s\sum_{l=0}^{N-1}\frac{1}{(jN+l)!}
\sum_{k=l}^\infty p_{jN+k}(jN+k)(jN+k-1)...(k-l+1) \nonumber \\ 
& & \times\sum_{t=1}^\infty\sum { ^\prime}
\prod_{u=1}^{jN+l}\tau_{N,n}(n_u)\prod_{v=jN+l+1}^{jN+k}
\gamma_{N,n}(n_v)s^{t-1}. \nonumber
\end{eqnarray}
 Observe that the coefficient of $s^{t-1}$ in the series
\[
\sum_{t=1}^\infty\sum { ^\prime}
\prod_{u=1}^{jN+l}\tau_{N,n}(n_u)\prod_{v=jN+l+1}^{jN+k}
\gamma_{N,n}(n_v)s^{t-1}
\]
can be written as
$$
\sum_{h=0}^{t-1}\ \sum_{n_1+...+n_{jN+l}=h}\
\prod_{u=1}^{jN+l}\tau_{N,n}(n_u)\
\sum_{n_{jN+l+1}+...+n_{jN+k}=t-1-h}\ \prod_{v=jN+l+1}^{jN+k}
\gamma_{N,n}(n_v).
$$
The rule of multiplying power series implies that this coefficient
equals the coefficient of $s^{t-1}$ in the power series expansion
of
\[\left[\sum_{i=1}^\infty\tau_{N,n}(i)s^i\right]^{jN+l}\left[\sum_{i=1}^\infty
\gamma_{N,n}(i)s^i\right]^{k-l} = [\psi_{N,n}(s)]^{jN+l}
[\phi_{N,n}(s)]^{k-l},
\]
where $\psi_{N,n}(s)$ and $\phi_{N,n}(s)$ are defined in
(\ref{psiphi}). Therefore,
\begin{eqnarray}
E(s^{\nu_{n+1}}; V_{N,n+1}=j) & = & s\sum_{l=0}^{N-1}
\frac{[\psi_{N,n}(s)]^{jN+l}}{(jN+l)!} \sum_{k=l}^\infty
p_{jN+k}(jN+k)(jN+k-1)...(k-l+1)[\phi_{N,n}(s)]^{k-l} \nonumber \\
 &= & s\sum_{l=0}^{N-1}\frac{[\psi_{N,n}(s)]^{jN+l}}{(jN+l)!}
f^{(jN+l)}(\phi_{N,n}(s)), \nonumber
\end{eqnarray}
which coincides with the right-hand side of (\ref{recurrence}).
This completes the proof.

\section{Fractional linear offspring}

Let $f(s)$ be a fractional linear pgf given by \be \label{flpgf}
f(s)=1-\frac{\DS b}{\DS 1-p}+\frac{\DS bs}{\DS 1-ps} \ee and the
parameter space $\{(p,b): 0<p<1, 0<b\leq 1-p\}$. Then the
offspring distribution is given by the geometric series
$p_k=bp^{k-1}, \ k=1,2,\ldots; p_0=1-\sum_{k=1}^\infty p_k$ and
the offspring mean is $m=b/(1-p)^2$. In the particular case
$b=p(1-p)$ we have $p_k=(1-p)p^k, \ k\geq 0$ which is the standard
geometric distribution with pgf $f(s)=(1-p)/(1-ps)$. It can be
verified, see Pakes and Dekking (1991), p. 361 if $N\geq 2$ and
Harris (1963), p. 9 if $N=1$, that for $N\in \cal N$ \be
\label{geotau} 1-p(1-\tau_N)=[b/(1-p)]^{1/N}[p\tau_N]^{1-1/N}. \ee

\vspace{0.5cm}{\bf Proposition 1}\ If the offspring distribution
has the fractional linear pgf (\ref{flpgf}), then $V_N$ follows a
zero-modified geometric (i.e., fractional linear) distribution
given by \be \label{flprobN}
P(V_N=j)=\frac{b}{p(1-p)}(1-\theta_N)\theta_N^j \quad (j\geq 1),
\quad P(V_N=0)=1-\frac{b}{p(1-p)}\theta_N\ee and \be
\label{flmean} EV_N=
\frac{b}{p(1-p)}\frac{\theta_N}{1-\theta_N},\ee where
\[
\theta_N=\left(\frac{p\tau_N}{1-p(1-\tau_N)}\right)^{N}
\]
and
 $\tau_N$ is the largest solution in $[0,1]$ of (\ref{geotau}).

\vspace{0.5cm}{\bf Proof}\ \ Since
$f^{(i)}(s)=i!bp^{i-1}/(1-ps)^{i+1} \ \ (i\geq 1)$, we have from
(\ref{main}) for $j\geq 1$ \beq P(V_N=j)
    & = &
    \sum_{k=0}^{N-1}\frac{\tau_N^{jN+k}}{(jN+k)!}\frac{b(jN+k)!p^{jN+k-1}}{(1-p(1-\tau_N))^{jN+k+1}} \nonumber \\
    & = &
    \frac{bp^{jN-1}\tau_N^{jN}}{(1-p(1-\tau_N))^{jN+1}}\sum_{k=0}^{N-1}\frac{(p\tau_N)^k}{(1-p(1-\tau_N))^k}.
    \nonumber
\eeq  Now, setting $(\theta_N)^{1/N}=p\tau_N/(1-p(1-\tau_N))$ one
can obtain the first formula in (\ref{flprobN}), which in turn
leads to (\ref{flprobN}) and (\ref{flmean}).

\vspace{0.5cm}{\bf Corollary}\ \  If the offspring distribution is
geometric, i.e., $p_k=(1-p)p^k, \ k\geq 0$, then  $V_N$ is
geometric as well,
$P(V_N=j)=(1-\tau_N)\tau_N^j \ (j\geq 0) \ \mbox{and} \
EV_N=\tau_N(1-\tau_N)^{-1}$, where $\tau_N$ is the largest solution
in $[0,1]$ of $ (\tau_N+1/m)^N=\tau_N^{N-1}\  \ (N\geq 1)$.

\vspace{0.5cm}{\bf Proof}\ \ In the case of geometric offspring
(\ref{flpgf}) holds with $b=p(1-p)$ and $m=p/(1-p)$. The equation
for $\tau_N$ follows by inspection from (\ref{geotau}). It is also
given in Pakes and Dekking (1991), p.361 if $N\geq2$. Simple
algebraic manipulations show that this equation simplifies to
$\theta_N=\tau_N$. Now, the rest of the statement follows from
(\ref{flprobN}) and (\ref{flmean}).

\vspace{0.5cm}{\it Remark}\ \ For geometric offspring with mean
$m>1$ we have $ P(V_1=j)=(1/m)(1-1/m)^j$ and $EV_1=m-1$. In
particular, $P(V_1=0)=1/m$ which equals the probability of
extinction, see
 Harris (1963), p. 9.

\vspace{0.5cm}Table~1 lists the probabilities $P(V_N=j)$,  $j=0,
1,2,\ldots 9$ as well as $EV_N$ for $1\leq N\leq 5$. The critical
mean values (see Section~1) are as follows: $m^c_1=1$, $m^c_2=4$,
$m^c_3=6.75$, $m^c_4=9.481$, $m^c_5=12.207$. The expected values
in the last column provide a measure of how many $N$-ary subtrees
$(1\leq N \leq 5)$ are supported by the geometric family tree with
offspring mean fixed to be $m=13$. See also Table~2 below for a
comparison with the Poisson offspring case.

\begin{table}[h]
\begin{center}
\begin{tabular}{|c|c|c|c|c|c|c|c|c|c|c|c|c|}
  \hline
  $V_N=$& 0 & 1& 2& 3& 4& 5& 6& 7& 8& 9& $\geq$10 & $E(V_N)$\\
  \hline
  $N=1$ & 0.08 & 0.07 & 0.07 & 0.06 & 0.06 & 0.05 & 0.05 & 0.04 & 0.04 & 0.04 & 0.44 & 12\\
  \hline
  $N=2$ & 0.16 & 0.14 & 0.11 & 0.10 & 0.08 & 0.07 & 0.06 & 0.05 & 0.04 & 0.03 & 0.16 & 5.22\\
  \hline
  $N=3$ & 0.26 & 0.19 & 0.14 & 0.11 & 0.08 & 0.06 & 0.04 & 0.03& 0.02& 0.02& 0.05 & 2.91 \\
  \hline
  $N=4$ & 0.37 & 0.23 & 0.15 & 0.09 & 0.06 & 0.04 & 0.02 & 0.01& 0.01& 0.01& 0.01 & 1.71\\
 \hline
  $N=5$ & 0.53 & 0.25 & 0.12 & 0.05 & 0.03 & 0.01 & 0.01  & 0& 0& 0& 0& 0.87 \\
  \hline
\end{tabular}
\end{center}
\caption{Probability distribution of $V_N$ assuming geometric
offspring with $m=13$.} \label{tabgeo}
\end{table}

\section{Poisson offspring} Consider the case of Poisson
offspring distribution with pgf given by \begin{equation}
\label{popgf} f(s)=e^{m(s-1)} \quad (m>0).
\end{equation}
 Then, the probability $\tau_N$ is the largest
solution of
\begin{equation} \label{poeq}
(1-s)e^{ms}=\sum_{j=0}^{N-1}(ms)^j/j!
\end{equation}
(see Pakes and Dekking (1991), p. 364). Since $f^{(i)}(s)=m^i
e^{m(s-1)} (i\ge 0)$, formula (\ref{main}) becomes
$$
\label{PoprobN} P(V_N=j)=e^{-m\tau_N}\sum_{k=0}^{N-1}
\frac{(m\tau_N)^{jN+k}}{(jN+k)!}, j\ge 0.
$$
Therefore, we have the following

{\bf Proposition 2}\ If the offspring distribution has the Poisson
pgf (\ref{popgf}), then
$$
P(V_N=j)=P(jN\le Y_{N}\le jN+N-1),
$$
where $Y_{N}$ has the Poisson pmf
$$
P(Y_{N}=k)=(m\tau_N)^k e^{-m\tau_N}/k! \qquad k=0,1,2,\ldots
$$
and $\tau_N$ is the largest solution in $[0,1]$ of equation
(\ref{poeq}).

Notice that $V_1$ has a Poisson distribution with parameter
$m\tau_1$. To calculate the critical value $m^c_N$ that yields a
non-zero solution $\tau^c_N$ in $[0,1]$ of equation (\ref{poeq})
we first notice that the product $y=m^c_N\tau^c_{N}$ satisfies the
equations
\begin{equation} \label{yeq2}
y^N/(N-1)!+\sum_{j=0}^{N-1}y^j/j!=e^y;
\end{equation}
 see Pakes and Dekking (1991), p. 365. Following their way of
 calculation, one can find $m_N^c$ and $\tau_N^c$ by substituting
 the solution of (\ref{yeq2}) into
\begin{equation} \label{yeq1}
my^{N-1}/(N-1)!=e^y.
\end{equation}
In case of binary trees, one can also use the Cayley's tree
function $y(z)=\sum_{k=1}^\infty k^{k-1}z^k/k!$ (see e.g. Odlyzko
(1995), Section 6.2) evaluated at $z=1/m_N^c$ for the solution of
(\ref{yeq2}). Inserting it into (\ref{yeq2}), we obtain
$m^c_{2}=3.3509$ and $\tau^c_{2}=0.5352$.

Our final remark concerns the case $m\to\infty$. It is easily seen
that Proposition~2 and the normal approximation of the Poisson
distribution imply a local limit theorem for $V_N$. Moreover,
Pakes and Dekking (1991) showed that in this case $\tau_N\to 1$.
This enables one to centralize and scale the limiting variable
$V_N$ in terms of the single parameter $m$ only.

Table~2 gives the probabilities $P(V_N=j)$,  $j=0, 1,2,\ldots 9$
as well as $EV_N$ for $2\leq N\leq 5$. The critical mean values
are as follows: $m^c_2=3.3509$, $m^c_3=5.1494$, $m^c_4=6.7993$,
$m^c_5=8.3653$.

\begin{table}[h]
\begin{center}
\begin{tabular}{|c|c|c|c|c|c|c|c|c|c|c|c|c|}
  \hline
  $V_N=$& 0 & 1& 2& 3& 4& 5& 6& 7& 8& 9& $\geq$10 & $E(V_N)$\\
  \hline
  $N=2$ & 0 & 0 & 0.01 & 0.04 & 0.11 & 0.19 & 0.22 & 0.19 & 0.13 & 0.07 & 0.04 & 6.25\\
  \hline
  $N=3$ & 0 & 0.01 & 0.09 & 0.25 & 0.32 & 0.22 & 0.08 & 0.02& 0& 0& 0.01& 4.00 \\
  \hline
  $N=4$ & 0 & 0.05 & 0.30 & 0.41 & 0.19 & 0.04 & 0 & 0& 0 & 0& 0.01& 2.87\\
 \hline
  $N=5$ & 0 & 0.17 & 0.51 & 0.28 & 0.04 & 0 & 0  & 0& 0& 0& 0& 2.19 \\
  \hline
\end{tabular}
\end{center}
\caption{Probability distribution of $V_N$ assuming Poisson
offspring with $m=13$.} \label{tabpoa}
\end{table}

\section{One-or-many offspring}

In this section we consider a two-parameter family of 1-or-r
offspring distributions defined for some $p\in (0,1)$ by $p_1=1-p$
and $p_r=p$, where $r>N>1$. Its pgf is $ f(s)=(1-p)s+ps^r $ and
thus $f'(s)=1-p+prs^{r-1}$ and $f^{(k)}(s)=pr(r-1)\ldots
(r-k+1)s^{r-k} \ (2\leq k\leq r)$. The probability $\tau_N$ is the
largest solution in $[0,1]$ of
\begin{equation} \label{eq16}
s=p\sum_{k=N}^r {r \choose k}s^k(1-s)^{r-k}
\end{equation}
(see again Pakes and Dekking (1991), p.366). Applying (\ref{main})
it is not difficult to obtain
\[
P(V_N=0)=1-p+p\sum_{k=0}^{N-1}{r \choose
k}\tau_N^k(1-\tau_N)^{r-k}
\]
and for $j=1,2,\ldots$ and $r\geq jN$
\[
P(V_N=j)=p\sum_{k=jN}^{jN+U}{r \choose k}\tau_N^k(1-\tau_N)^{r-k},
\]
where $U=\min\{N-1, \ r-jN\}$. Let $B_r(\tau_N)$ denote a binomial
$(r, \tau_N)$ random variable.

{\bf Proposition 3}\ If the offspring pgf is $f(s)=(1-p)s+ps^r $
$(1\leq N <r)$ and $\tau_N$ is the largest solution in $[0,1]$ of
(\ref{eq16}), then $P(V_N=0)=1-p+pP(B_r(\tau_N)\leq N-1)$ and for
$j=1,2,\ldots$
 \be \label{omprobN}
P(V_N=j)= pP(jN \leq B_r(\tau_N)\leq jN+U) \quad \mbox{if} \
jN\leq r, \ee where $U=\min\{ N-1, \ r-jN\}$ and $P(V_N=j)= 0$ if
$jN>r$. The expected value of $V_N$ is
\[
EV_N=p\sum_{j=1}^{[r/N]}jP(jN \leq B_r(\tau_N)\leq jN+U),
\]
where $[x]$ is the integer part of $x$.

In particular, if $r=N+1$ or $r=N+2$ and $N>2$, then
(\ref{omprobN}) implies that $V_N$ takes on values 0 or 1; if
$N=2$ and $r=4$, then $V_N$ takes on values 0, 1, or 2. Table~3
provides some numerical illustrations. Note that the offspring
mean $m=13.09$ enables comparisons with Tables 1 and 2.

\begin{table}[h]
\begin{center}
\begin{tabular}{|c|c|c|c|c|c|c|c|c|c|}
  \hline
  $V_N=$& 0 & 1& 2& 3& 4& 5& 6& 7& $E(V_N)$\\
  \hline
  $N=2$ & 0.07 & 0 & 0& 0 & 0& 0.06 & 0.53 & 0.34 & 5.86 \\
  \hline
  $N=3$ & 0.07 & 0 & 0 & 0.06 & 0.87 & 0 & 0  & 0& 5.05 \\
  \hline
  $N=4$ & 0.07 & 0 & 0.06 & 0.87 & 0 & 0 & 0 & 0& 2.73 \\
 \hline
  $N=5$ & 0.07 & 0 & 0.93 & 0 & 0 & 0 & 0  & 0& 1.86 \\
  \hline
\end{tabular}
\end{center}
\caption{Probability distribution of $V_N$ assuming 1-or-14
offspring with $p=0.93$ ($m=13.09$).} \label{taboom}
\end{table}

It is interesting to point out the following relationship between
the 1-or-r and Poisson offspring cases. There exists (see Pakes
and Dekking (1991)) a critical value $p_N^c$ such that for
$p=p_N^c$ equation (\ref{eq16}) has a single solution $\tau_N^c$
in $(0,1)$. Suppose that $\lim_{r\to \infty}(r\tau_N^c)\to y$,
where $y$ satisfies (\ref{yeq1}) and (\ref{yeq2}). Then, applying
Theorem~7, Pakes and Dekking (1991), one can obtain that $V_N(r)$
converges in distribution to $V_N(y)$, where $V_N(r)$ and $V_N(y)$
are copies of $V_N$ assuming one-or-many and Poisson offspring
with mean $m_N^c$, respectively.

\section*{Acknowledgments}

We thank the referee for his valuable comments and suggestions and
especially for his help to eliminate some defects in Proposition
1. This work was done during L. Mutafchiev's visit at the
Mathematics Department of the University of South Florida in
2004-05 academic year. He thanks for the hospitality and support.
G. Yanev is partially supported by NFSI-Bulgaria, MM-1101/2001.

\end{document}